\documentclass[french,openright,12pt,twoside]{amsart}
\usepackage{a4}
\usepackage[T1]{fontenc}
\usepackage[all]{xy}
\usepackage{babel,indentfirst}
\usepackage{verbatim}
\usepackage{stmaryrd}
\usepackage{amsfonts}
\usepackage{amsmath,amscd,amssymb}
\usepackage{amsthm}
\usepackage{enumerate}
\usepackage{xspace}
\usepackage{ulem}
\usepackage{euscript}
\usepackage{epsfig}

\newtheorem {prop}{Proposition}

\newcommand  {\GL }{\mathop{\mathrm{GL}}\nolimits}

\newcommand {\Vect}{\mathop{\mathrm{Vect}}\nolimits}

\newcommand  {\card }{\mathop{\mathrm{Card}}\nolimits}

\newcommand  {\qu}{\mbox{ }}
\newcommand {\V}{\par\vspace {6pt}}
\newcommand {\tir}{\quad - \mbox{ }}

\begin {document}
\title{Une propriété du groupe à $168$ éléments}
\author{R. Goblot}
\date{24 mars 2003}

{\thanks {Classification AMS: Arithmetic and combinatorial problem 20 D 
60, Simple groups 20 E 32, Geometric grouptheory 20F32}}
\maketitle 

Soit $K$ un corps commutatif, 
$\mathcal E$ un espace affine sur $K$ de dimension $n$, $E$ l'espace des 
vecteurs de $\mathcal E$, $\mathcal G$ le groupe affine et $\mathcal T$ 
le sous-groupe des translations de $\mathcal E$, 
$G\sim \mathcal G/\mathcal T$ le groupe linéaire de $E$. 
Pour tout $\varphi \in \mathcal G$, on 
note $\overline \varphi $ l'application linéaire associée dans $G$. Pour 
tout $f\in G$ et tout $a\in \mathcal E$, on note $f_a$ l'application 
affine fixant $a$ d'application linéaire associée $f$. À tout $a\in 
\mathcal E$ correspond donc la section 
$$\mathrm s_a\colon G\hookrightarrow \mathcal G\qu,\qu f\mapsto f_a$$
L'objet de ce travail est de montrer le théorème suivant:\V

\noindent\textbf{Théorème} \quad 
$a)$ \textit{Toute section $s\colon G\rightarrow \mathcal G$ est du type 
$\mathrm s_a$ à l'exception du cas où $K=\mathbf F_2$ et $n=3$.}

$b)$ \textit{On suppose $K=\mathbf F_2$ et $n=3$. Soit $\mathcal F$ 
l'ensemble des sous-groupes de $G$ de cardinal $7$. Pour tout $M\in 
\mathcal F$, soit $\mathrm N(M)$ le normalisateur de $M$. L'ensemble 
$\mathcal F$ est muni d'une structure canonique de $\mathbf F_2$-espace 
affine de dimension $3$ d'espace des vecteurs $E$: pour $M_1$ et $M_2$ 
distincts de 
$\mathcal F$, $\mathrm N(M_1)\cap \mathrm N(M_2)=\{\mathrm I,g,g^2\}$
 est un sous-groupe de $G$ d'ordre $3$. On pose $\overrightarrow 
 {M_1M_2}=\overrightarrow \varepsilon $ où $\overrightarrow \varepsilon $ 
 est le vecteur de $E$ fixe par $g$. L'opération de $G$ sur $\mathcal F$ 
 par conjugaison donne une section de $G$ dans le groupe affine de 
 $\mathcal F$ sans point fixe.}
 
 $c)$ \textit{L'ensemble $ \mathcal E'$ des isomorphismes affines 
 $\mathcal E \rightarrow \mathcal F$ induisant l'identité sur $E$ a une 
 structure canonique de $\mathbf F_2$-espace affine de dimension $3$ 
 de même espace vectoriel associé $E$ et de même groupe affine $\mathcal G$ que 
 $\mathcal E$. 
 L'ensemble des sections $G\hookrightarrow \mathcal G$, de cardinal 
 $16$, est en bijection canonique avec $\mathcal E\cup \mathcal E'$. Les 
 espaces affines $\mathcal E$ et $\mathcal E'$ jouent le même rôle l'un 
 par rapport à l'autre.}

\paragraph{Cas où $K$ n'est pas $\mathbf F_2$}
Il s'agit de montrer \textit{si $K\ne \mathbf F_2$, toute section est du 
type $\mathrm s_a$ où $a\in 
\mathcal E$.} 

Soit 
$s $ une section et $\lambda \ne 1$ dans $K^{\times}$. L'image par 
$s $ de l'homothétie vectorielle $\lambda \mathrm I_{E}$ se remonte 
par $s $ en une homothétie affine  $\mathrm{hom}(a,\lambda )$ de 
rapport $\lambda $ et de centre $a $. Soit  $f\in G$ et $\varphi =s(f)$.
 L'homothétie vectorielle $\lambda 
\mathrm I_E$ commute avec $f$, donc $s (\lambda \mathrm 
I_E)=\mathrm{hom}(a,\lambda ) $ 
commute avec $s(f)=\varphi $. On a donc 
$$\varphi (a )= \bigl(\varphi \circ \mathrm{hom}(a,\lambda )\bigr)(a )= 
\bigl(\mathrm{hom}(a,\lambda )\circ \varphi \bigr)(a)=
\bigl(\mathrm{hom}(a,\lambda )\bigr)\bigl(\varphi (a )\bigr)$$
donc l'homothétie affine $\mathrm{hom} (a,\lambda )$ laisse fixes
$a $ et $\varphi (a )$. On a donc $\varphi (a)=a 
$ et $s $ est la section $\mathrm s_a$ envoyant $G$ en le stabilisateur 
de $a $. Ce raisonnement tombe en défaut si $K=\mathbf F_2$ car il n'existe 
aucun $\lambda \ne 1$.\V 

Dans toute la suite, le corps est supposé être $\mathbf F_2$. Dans la 
section 1, on étudie les sous-groupes de $G$. On sait que le groupe à 
$168$ éléments se présente aussi sous l'aspect de $\mathrm{PSL}(2,\mathbf 
F_7)$. Ceci est évoqué sans démonstration dans la remarque 1.7. Les 
sections 2 et 3 développent respectivement les points $b)$ et $c)$ du 
théorème. La section 4 étudie les $\mathbf F_2$-espaces affines de 
dimension supérieure ou égale à $4$.

\section{Sous groupes du groupe $G$}

\paragraph{1.1 \qu Plans projectifs sur $\mathbf F_2$}
Le cardinal de $G$ est celui de l'ensemble des bases de $E$, soit 
$7\times 6\times 4=168$.
Le groupe multiplicatif $\mathbf F_2^{\times}$ étant réduit à l'unité, la géométrie 
vectorielle en dimension $3$ sur $ E$ 
et la géométrie projective en dimension $2$ sur le plan projectif 
$\mathrm P( E)$ 
se confondent et on a une bijection canonique
$$E\setminus \{\overrightarrow 0\}\rightarrow \mathrm P(E)\qu,\qu 
\overrightarrow v\mapsto v$$
On notera de la même façon une application linéaire inversible
 $f\colon E\rightarrow E$ et 
l'homographie $f\colon \mathrm P(E)\rightarrow \mathrm P(E)$ associée. 
La  donnée d'une base $(\overrightarrow {e_1},\overrightarrow 
{e_2},\overrightarrow {e_3})$ de $E$ équivaut à la donnée d'un triangle 
non aplati $(e_1,e_2,e_3)$ de $\mathrm P(E)$, qui se complète d'unique 
façon en un repère $(e_1,e_2,e_3,e_4)$ de $\mathrm P(E)$ où $e_4$ est 
l'unique point de $\mathrm P(E)$ ne se trouvant sur aucun des côtés du 
triangle. Une homographie est donc définie par la donnée d'un triangle 
et du triangle image.

Le plan projectif $\mathrm P(E)$ comporte $7$ points et $7$ droites. 
Par chaque point passent $3$ droites et sur chaque droite on a $3$ 
points. La figure suivante constituée d'un triangle $abc$ et de $3$
droites  $(ap),(bq),(cr)$ concourantes en $o$ permet de se 
représenter un plan projectif sur $\mathbf F_2$ en convenant de 
considérer $p,q,r$ comme alignés. L'ensemble des $7$ droites de $\mathrm P(E)$
constitue le plan projectif dual $\mathrm P(E^*)=\mathrm P^*(E)$. 

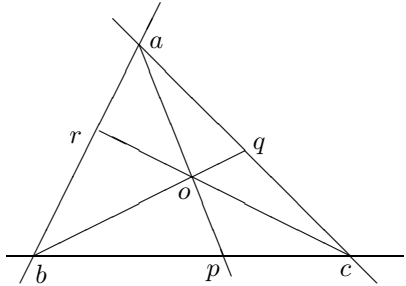
\begin{center}
\begin{picture}(140,100)
\put (0,10){\line(1,0){150}}
\put (5,0){\line(1,2){53}}
\put (40,100){\line(1,-1){100}}

\put (10,10){\line(2,1){80}}
\put (130,10){\line(-2,1){95}}
\put (50,90){\line(2,-5){35}}

\put(54,91){\makebox(0,0)[l]{{\footnotesize $a$}}}
\put(13,7){\makebox(0,0)[t]{{\footnotesize $b$}}}
\put(128,7){\makebox(0,0)[t]{{\footnotesize $c$}}}

\put(78,7){\makebox(0,0)[t]{{\footnotesize $p$}}}
\put(93,52){\makebox(0,0)[l]{{\footnotesize $q$}}}
\put(26,57){\makebox(0,0)[t]{{\footnotesize $r$}}}
\put(67,36){\makebox(0,0)[t]{{\footnotesize $o $}}}

\end{picture}
\end{center}

\begin{prop} Soit $P$ un ensemble de  \textit{points}, 
$P^*$ un ensemble de cardinal supérieur à $1$ de
\textit{droites}.
 On suppose vérifiées les conditions d'incidence suivantes:

\tir pour toute  droite, il existe exactement trois points en incidence avec cette 
droite,

\tir pour toute paire de deux points distincts, il existe   une unique droite en 
incidence avec ces deux points, 

\tir pour toute paire de deux droites distinctes, il existe   un unique point en 
incidence avec ces deux droites.

\noindent Alors
 $P$ a une structure de $\mathbf F_2$-plan projectif et 
$P^*$ s'identifie à l'ensemble des droites de $P$.
Pour qu'une bijection $h\colon P\rightarrow P$ soit une homographie, il 
faut et il suffit qu'elle conserve l'alignement.\end{prop}

Montrons qu'un point $o$ est en incidence avec exactement $3$ droites. 
Soit $m_1$ 
distinct de $o$. Comme $(om_1)$ n'est pas la seule droite, il existe un 
point $m_2$ non sur $(om_1)$. Soit $m_3$ le troisème point de la droite 
$D=(m_1m_2)$. 
Les trois droites $(om_1)$, $(om_2)$, $(om_3)$
 sont en incidence avec $o$. 
S'il existait une quatrième droite, elle couperait $D$ en un quatrième 
point ce qui est impossible. Pour tout $i\in \{1,2,3\}$, soit $m'_i$ le 
troisème point de la droite $(om_i)$, alors 
$P=\{o,m_1,m'_1,m_2,m'_2,m_3,m'_3\}$ est de cardinal $7$.

Adjoignons un élément $\overrightarrow 0$ à $P$.  Un élément sera noté $x$ 
ou $\overrightarrow x$ selon qu'on le considère comme appartenant à $P$ 
ou à $E=P\cup \{\overrightarrow 0\}$. On définit sur $E$ l'addition suivante:

\tir $\forall \overrightarrow x\in E\qu, \qu \overrightarrow 
x+\overrightarrow 0=\overrightarrow 0+\overrightarrow x=\overrightarrow x$,

\tir $\forall \overrightarrow x\in E\qu, \qu \overrightarrow 
x+\overrightarrow x=\overrightarrow 0$,

\tir pour $\overrightarrow x$ et $\overrightarrow y$ distincts de 
$\overrightarrow 0$ et entre eux, $\overrightarrow x+\overrightarrow 
y=\overrightarrow z$ où $z$ est le troisième point de la droite en 
incidence avec $x$ et $y$.

\noindent Montrons que cette addition est associative: $(\overrightarrow \alpha 
+\overrightarrow \beta )+\overrightarrow \gamma =\overrightarrow \alpha 
+(\overrightarrow \beta +\overrightarrow \gamma )$. Regardons le cas où 
$\alpha ,\beta ,\gamma $ sont distincts non alignés. On pose 
$\overrightarrow \delta =\overrightarrow \alpha +\overrightarrow \beta $ 
et $\overrightarrow \varepsilon =\overrightarrow \beta +\overrightarrow 
\gamma $. Alors 
*\begin{eqnarray} (\overrightarrow \alpha +\overrightarrow \beta 
)+\overrightarrow \gamma &=& \overrightarrow \delta +\overrightarrow 
\gamma =\overrightarrow \sigma \nonumber\\
\overrightarrow \alpha +(\overrightarrow \beta +\overrightarrow \gamma 
)&=&\overrightarrow \alpha +\overrightarrow \varepsilon =\overrightarrow 
\sigma \nonumber\end{eqnarray}
où $\sigma $ est l'unique point en incidence avec  les
 droites  $(\delta \gamma )$ et $(\alpha \varepsilon )$.

Alors $E$, de cardinal $8$,
 a une structure de $\mathbf F_2$-espace vectoriel de dimension 
$3$, $P$ s'identifie au plan $\mathrm P(E)$ et $P^*$ au dual de $P$.

Montrons que si une bijection $h\colon P\rightarrow P$ conserve 
l'alignement, alors c'est une homographie. On prolonge $h$ à $E$ en 
posant $h(\overrightarrow 0)=\overrightarrow 0$. Soit $\overrightarrow x$,
 $\overrightarrow y$  non nuls et distincts dans $E$, 
 $\overrightarrow z=\overrightarrow x+\overrightarrow y$. Alors $x,y,z$ 
 sont alignés distincts, donc $h(x),h(y),h(z)$ sont alignés distincts, d'où 
 $h(\overrightarrow x)+h(\overrightarrow y)=h(\overrightarrow z)$ et 
 $h\colon E\rightarrow E$ est bien linéaire.

\paragraph{1.2\qu Sous-groupes isomorphes à $\mathbf S_4$} Soit $D$ une droite 
de $\mathrm P(E)$. L'ensemble $\mathrm P(E)\setminus D$ est formé de 
quatre points et constitue un repère projectif de $\mathrm P(E)$. 
Considérant l'action du stabilisateur $\mathrm{St}(D)$ de $D$ sur 
le quadrilatère $\mathrm P(E)\setminus D$, on obtient un isomorphisme
$\mathrm{St}(D)\sim \mathbf S_4$. Soit $G(D)\sim \mathbf S_3$ le groupe 
des homographies de $D$. L'application \textquotedblleft{restriction}\textquotedblright\ à $D$ définit 
un morphisme $\mathrm{St}(D)\sim \mathbf S_4\rightarrow G(D)\sim \mathbf 
S_3$ dont le noyau est un groupe de Klein $\mathrm{Kl}(D)$: l'ensemble des $h\in G$ 
induisant l'identité sur $D$. Ce groupe $\mathrm{Kl}(D)$  est formé de l'identité 
et de trois involutions, chacune induisant deux transpositions 
disjointes sur le quadrilatère  $\mathrm P(E)\setminus D$. On a $7$ droites 
dans $\mathrm P(E)$, une classe de conjugaison de $7$ groupes 
$\mathrm{St}(D)$ et une classe de conjugaison de $7$ groupes de Klein 
$\mathrm{Kl}(D)$.

Dualement, le stabilisateur $\mathrm{St}(\varepsilon )$ 
d'un point $\varepsilon $ est isomorphe à 
$\mathbf S _4$. Le sous-groupe de Klein $\mathrm{Kl}(\varepsilon )$ est formé des 
homographies laissant stable chacune des trois droites passant par 
$\varepsilon $. Une involution de $\mathrm{Kl}(\varepsilon )$ laisse fixe 
les points d'une des 
trois droites passant par $\varepsilon $ et  échange les deux 
points distincts de $\varepsilon $ sur chacune des deux autres droites.
On a $7$ points 
dans $\mathrm P(E)$, une classe de conjugaison de $7$ groupes 
$\mathrm{St}(\varepsilon )$ et une classe de conjugaison de $7$ groupes de Klein 
$\mathrm{Kl}(\varepsilon )$.

\paragraph{1.3\qu Sous-groupes d'ordre $8$ et involutions} Il s'agit de la classe de 
conjugaison des  $2$-Sylow(s). Dans chaque groupe isomorphe à $\mathbf 
S_4$, on a trois sous-groupes diédraux de cardinal $8$. Un $2$-Sylow est 
donc contenu dans $\mathrm{St}(D)\cap \mathrm{St}(\varepsilon )$  où $D$ 
est une droite et $\varepsilon $ un point. Nécessairement, $\varepsilon 
\in D$, sinon $\varepsilon $ est un point du repère $\mathrm 
P(E)\setminus D$ et $\mathrm{St}(D)\cap \mathrm{St}(\varepsilon )\sim 
\mathbf S_3$ est le 
stabilisateur de $\varepsilon $ dans $\mathrm{St}(D)\sim \mathbf S_4$. 

\textit{Un $2$-Sylow est donc le stabilisateur $\mathrm{St}(D,\varepsilon )$
d'un couple $(D,\varepsilon )$ 
où $D$ est une droite et $\varepsilon $ un point de $D$.} On a $7$ droites 
et sur chacune $3$ points. On a donc $21$ couples $(D,\varepsilon )$. La 
classe de conjugaison des $2$-Sylow(s) est de cardinal $21$. 

Toute involution étant contenue dans un $2$-Sylow, les involutions de 
$G$ sont celles  mises en évidence précédemment: ce sont des 
paires de transpositions disjointes $(a,b),(c,d)$ où le quadriltaère 
$a,b,c,d$ est un repère projectif de $\mathrm P(E)$. On a $7$ tels 
quadrilatères et pour chacun d'eux $3$ telles involutions. Les 
involutions de $G$ forment une classe de conjugaison de cardinal $21$.

On a une bijection canonique de l'ensemble des $21$ $2$-Sylow(s) sur 
l'ensemble des $21$ involutions: à tout $2$-Sylow 
$\mathrm{St}(D,\varepsilon )$, on associe l'involution 
$\tau _{D,\varepsilon }$ telle que $\{\mathrm I ,\tau_{D,\varepsilon }\}$ 
soit le centre de $\mathrm{St}(D,\varepsilon )$. Si $a,b,c,d$ sont les 
$4$ points de $\mathrm P(E)\setminus D$, choisissant les notations de 
sorte que les droites $(ab)$ et $(cd)$ se coupent en $\varepsilon $, 
alors $\tau _{D,\varepsilon }$ échange $a$ et $b$ et échange $c$ et $d$.
On a $\{\mathrm I,\tau _{D,\varepsilon }\}=
 \mathrm{Kl}(D)\cap \mathrm{Kl}(\varepsilon )$.

\paragraph{1.4\qu Sous-groupes d'ordre $3$ et $6$} Soit $D$ une droite de 
$\mathrm P(E)$. Alors $\mathrm P(E)\setminus D$ est un quadrilatère et le 
groupe $\mathrm{St}(D)$ est isomorphe à $\mathbf S_4$. À tout 
$\varepsilon \in \mathrm P(E)\setminus D$ est associé le groupe 
$\mathrm{St}(D,\varepsilon )\sim \mathbf S_3$ laissant stables $D$ et 
$\varepsilon$ et le sous-groupe $\Gamma (D,\varepsilon )\sim \mathbf A_3$  
formé de deux permutations circulaires sur les $3$ points de $D$ et sur  
les $3$ points du triangle complémentaire de $D\cup\{\varepsilon \}$. Les groupes 
$\Gamma(D,\varepsilon )$ sont des $3$-Sylow(s). Comme les $3$-Sylow(s) 
sont conjugués, ils sont tous de ce type. 

On a $7$ possibilités pour la 
droite $D$ et $7-3=4$ possibilités pour le point $\varepsilon \notin D$. 
On a donc $7\times 4=28$ sous-groupes du type $\mathrm{St}(D,\varepsilon 
)\sim \mathbf S_3$ et $28$ sous-groupes $\Gamma(D,\varepsilon 
)\sim \mathbf A_3$. 

\paragraph{1.5\qu  Les classes de conjugaison des groupes isomorphes à 
$\mathbf S_4$ comme $\mathbf F_2$-plans projectifs en dualité} Les 
sous-groupes de $G$ isomorphes à $\mathbf S_4$ se rangent en deux classes 
de conjugaison: les $7$ stabilisateurs de point $\mathrm{St}(\varepsilon 
)$ où $\varepsilon \in P$, les $7$ stabilisateurs de droite 
$\mathrm{St}(D)$ où $D\in P^*$. On trouve les relations d'incidence:

 \tir Soit $A,B$ deux droites distinctes se coupant en 
 $\varepsilon $, $C$ la troisième droite passant par $\varepsilon $.
  On a $\mathrm{St}(A)\cap 
 \mathrm{St}(B)=\mathrm{St}(C,\varepsilon )=\{\mathrm I,\tau 
 _{C,\varepsilon }\}$. De même, 
 soit $a,b$ deux points distincts de $\mathrm P(E)$, $c$ le troisième 
 point de la droite $D=(ab)$. Alors $\mathrm{St}(a)\cap 
 \mathrm{St}(b)=\mathrm{St}(D,c)=\{\mathrm I,\tau _{D,c}\}$.
 
 \tir Soit $D$ une droite, $\varepsilon \notin D$ un point, alors 
 $\mathrm{St}(D)\cap \mathrm{St}(\varepsilon )=
 \mathrm{St}(D,\varepsilon )\sim \mathrm S_3$.

\tir Soit $D$ une droite, $\varepsilon\in D $ un point, alors 
$\mathrm{St}(D)\cap \mathrm{St}(\varepsilon )=\mathrm{St}(D,\varepsilon 
)$ est diédral de cardinal $8$.

\begin{prop} Soit $P,P^*$ les deux classes de conjugaison de sous-groupes 
isomorphes à $\mathbf S_4$. Soit $S,S'$ deux sous-groupes de $P\cup P^*$. 

\tir $S\cap S'$ est de cardinal $2$ si et seulement si
$S$ et $S'$ sont dans la même classe de conjugaison.

\tir $S\cap S'\sim \mathbf S_3$ si et seulement si $S$ et $S'$ ne sont 
pas dans la même classe de conjugaison et ne sont pas en incidence.

\tir $S\cap S'$ est diédral de cardinal $8$ 
si et seulement si $S$ et $S'$ ne sont 
pas dans la même classe de conjugaison et  sont  en incidence.\end{prop}

\paragraph{1.6\qu  Sous-groupes d'ordre $7$ de $G$} Ils constituent 
l'ensemble $\mathcal F$ des $7$-Sylow(s). 
On va raisonner dans le cadre vectoriel.
L'extension $\mathbf F_8$ de $\mathbf F_2$ est une $\mathbf F_2$-algèbre 
de rang $3$. Le groupe multiplicatif $\mathbf F_8^{\times}$ est 
cyclique d'ordre $7$. Les éléments de $\mathbf F_8^{\times}$ se rangent 
en deux classes $x,x^2,x^4$ et $x^{-1},x^{-2},x^{-4}$, racines des deux 
polynômes irréductibles de degré $3$: 
$$P(X)=X^3+X+1\quad ,\quad  Q(X)= X^3+X^2+1$$

Dans $G$, on a deux classes  $\mathcal P$ et 
$\mathcal Q$ d'éléments d'ordre $7$ 
selon que le polynôme minimal est $P(X)$ ou $Q(X)$. 

Soit $\mathcal B( E)$ l'ensemble des bases de $ E$. On a 
une application $\mathcal B( E)\rightarrow \mathcal P$ associant 
à la base $(\overrightarrow {e_1},\overrightarrow {e_2},\overrightarrow 
{e_3})$
 l'application dont la matrice dans cette base est la 
matrice compagnon de $P(X)$:
$$\left( \begin{array}{ccc} 0&0&1\\ 1&0&1\\ 0&1&0\end{array}\right)$$
Cette application est surjective. En effet, soit $f$ admettant $P(X)$ pour 
polynôme caractéristique. Comme $P(X)$ est irréductible, $f$ n'a aucun 
vecteur propre et admet la matrice compagnon de $P(X)$ dans toute base 
$\bigl (\overrightarrow {e_1},f(\overrightarrow {e_1}),
f^2(\overrightarrow {e_1})\bigr)$. 
Ainsi, $\mathcal P$ et $\mathcal Q$ 
sont les deux classes de conjugaison des éléments d'ordre $7$ de $G$.

Soit $\overrightarrow {e_1}\ne \overrightarrow 0$. L'application $f\mapsto 
\bigl(f(\overrightarrow {e_1}),f^2(\overrightarrow {e_1})\bigr)$ 
est donc une bijection de $\mathcal P$ sur 
l'ensemble des couples de vecteurs $(\overrightarrow 
{e_2},\overrightarrow {e_3})$ tels que
 $(\overrightarrow {e_1},\overrightarrow {e_2},\overrightarrow {e_3})$ 
soit base de $ E$. On a $8-2=6$ possibilités pour $\overrightarrow {e_2}$ et 
$8-4=4$ possibilités pour $\overrightarrow {e_3}$.
 Ainsi, $\mathrm{Card}(\mathcal P)= 
6\times 4=24$. Les éléments de $\mathcal P$ se rangent en $8$ paquets de 
$3$ élements $(f,f^2,f^{4})$ engendrant le même groupe $M=\langle 
f\rangle =\langle f^2\rangle= \langle f^4\rangle$ d'ordre $7$. L'ensemble 
$\mathcal F$ des sous-groupes d'ordre $7$ est donc de cardinal $8$.

 Considérons un tel sous-groupe $M$ et ses générateurs $f,f^2,f^4$ 
 appartenant à $\mathcal P$. Comme le coefficient de $X^2$ dans $P$ est nul,
  pour tout $\overrightarrow \varepsilon \ne \overrightarrow 0\in  E$,
   on a $f(\overrightarrow \varepsilon 
 )+f^2(\overrightarrow \varepsilon )+f^4(\overrightarrow \varepsilon )=
 \overrightarrow 0$. Considérant la 
 droite vectorielle 
 $\Delta=\{\overrightarrow 0,\overrightarrow \varepsilon\} $ et le plan vectoriel $\Pi= 
\bigl\{ \overrightarrow  0,f(\overrightarrow \varepsilon ),
f^2(\overrightarrow \varepsilon ),
f^4(\overrightarrow \varepsilon )\bigr\}$, on a $ E= \Delta \oplus \Pi$.
 La donnée du groupe 
$M=\langle f\rangle=\langle f^2\rangle =\langle f^4\rangle$ induit une 
orientation sur le plan $\Pi$ en prenant la permutation circulaire 
$\bigl(f(\overrightarrow \varepsilon ),f^2(\overrightarrow \varepsilon ),
f^4(\overrightarrow \varepsilon )\bigr)$, 
indépendante du générateur $f\in \mathcal P\cap M$. Ainsi, \textit{la 
donnée d'un $M\in \mathcal F$
détermine $7$ décompositions orientées $ E=\Delta \oplus \Pi$
où $\Delta =\{\overrightarrow 0,\overrightarrow \varepsilon \}$ et
 $\Pi=\bigl\{\overrightarrow 0,f(\overrightarrow \varepsilon 
),f^2(\overrightarrow \varepsilon ),f^4(\overrightarrow \varepsilon )\bigr\}$ muni de la permutation 
circulaire $\bigl (f( \overrightarrow \varepsilon ),
 f^2( \overrightarrow \varepsilon ),f^4(\overrightarrow  \varepsilon 
)\bigr)$ avec $f\in \mathcal P\cap M$.}

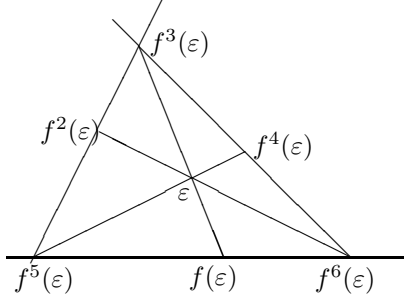
\begin{center}
\begin{picture}(140,100)
\put (0,10){\line(1,0){150}}
\put (9,8){\line(1,2){50}}
\put (40,100){\line(1,-1){90}}

\put (10,10){\line(2,1){80}}
\put (130,10){\line(-2,1){95}}
\put (50,90){\line(2,-5){32}}

\put(54,91){\makebox(0,0)[l]{{\footnotesize $f^3(\varepsilon )$}}}
\put(14,7){\makebox(0,0)[t]{{\footnotesize $f^5(\varepsilon )$}}}
\put(128,7){\makebox(0,0)[t]{{\footnotesize $f^6(\varepsilon )$}}}

\put(78,7){\makebox(0,0)[t]{{\footnotesize $f(\varepsilon )$}}}
\put(93,52){\makebox(0,0)[l]{{\footnotesize $f^4(\varepsilon )$}}}
\put(24,63){\makebox(0,0)[t]{{\footnotesize $f^2(\varepsilon )$}}}
\put(67,36){\makebox(0,0)[t]{{\footnotesize $\varepsilon  $}}}

\end{picture}
\end{center}
 
 \noindent Pour tout $k\in \mathbf Z$, les points $f^k(\varepsilon 
 ),f^{k+1}(\varepsilon ),f^{k+3}(\varepsilon )$ sont alignés.
 
 \paragraph{1.7\qu Remarques}\quad Le groupe $G$ donné abstraitement peut 
être considéré de deux façons comme groupe des homographies d'un $\mathbf 
F_2$-plan projectif $P$ selon que l'on choisit pour $P$ l'une ou l'autre 
des deux classes de conjugaison de sous-groupes isomorphes à $\mathbf 
S_4$. Le groupe  $\mathrm{Aut}(G)$ des automorphismes de $G$
 est de cardinal $336$, un 
automorphisme étant intérieur ou non selon qu'il laisse stables ou 
échange ces deux classes de conjugaison. 

Ceci permet de montrer que $G$ apparait aussi 
comme isomorphe à $\mathrm{PSL}(2,\mathbf F_7)$, $\mathrm{Aut}(G)$ étant 
alors isomorphe à $\mathrm{PGL}(2,\mathbf F_7)$. Dans $\mathbf F_7$, 
$j=2$ et $j^2=4$ sont racines cubique de l'unité. Soit $\Delta$ est une 
$\mathbf F_7$-droite projective. Appelons \textit{tétraèdre} tout 
quadruplet de birapport $-j$ ou $-j^2$, i.e. de stabilisateur isomorphe à 
$\mathbf A_4$. Si $[a_1,a_2,a_3,a_4]=-j$, on 
définit $a'_1,a'_2,a'_3,a'_4$ par les égalités
$$[a'_1,a_2,a_3,a_4]=[a_1,a'_2,a_3,a_4]=
[a_1,a_2,a'_3,a_4]=[a_1,a_2,a_3,a'_4]=-j^2$$
On a alors le \textit{tétraèdre symétrique} $[a'_1,a'_2,a'_3,a'_4]=-j^2$. 
L'adjonction des tétraèdres symétriques $A=(a_1,a_2,a_3,a_4)$
 et $A'=(a'_1,a'_2,a'_3,a'_4)$ forme le
\textit{cube} $(A,A')=(a_1,a_2,a_3,a_4)(a'_1,a'_2,a'_3,a'_4)$. 
On montre que l'on 
a $14$ cubes, i.e. $14$ façons de partager $\Delta$ en deux tétraèdres 
symétriques. L'action de $\mathrm{PSL}(\Delta)$ définit deux orbites 
$P$ et $P^*$ sur l'ensemble des cubes 
avec des relations d'incidence permettant de considérer 
$P,P^*$ comme deux $\mathbf F_2$-plans projectifs et 
$\mathrm{PSL}(\Delta)$ comme le groupe des homographies de $P$ et $P^*$.
 Par exemple, on montre que deux 
cubes $(A,A')$ et $(B,B')$ sont dans deux orbites distinctes et non 
en incidence si et seulement si $\card A\cap B=\card A'\cap B'$ et $\card 
A\cap B'=\card A'\cap B$ sont impairs, i.e. valent $1$ ou $3$.

Dans cette optique, on a une bijection canonique
$\Delta \rightarrow \mathcal F$, $m\mapsto M$ où le normalisateur $\mathrm N(M)$ de 
$M$ est le stabilisteur du point $m$ de $\Delta$. Cette bijection 
détermine sur $\mathcal F$ une structure canonique de $\mathbf 
F_7$-droite projective. L'action de $G$ sur $\mathcal F$ par conjugaison 
permet d'identifier $G$ à $\mathrm{PSL}(\mathcal F)$.

On a deux types d'actions fidèles de $G$ sur un ensemble $X$ de cardinal 
$8$:

\tir une action à deux orbites de cardinal $1$ et $7$ définissant sur $X$ 
une structure de $\mathbf F_2$-espace vectoriel de dimension $3$,

\tir une action transitive déterminant sur $X$ une structure de $\mathbf 
F_7$-droite projective par une bijection $\mathcal 
F\rightarrow X$.

\section{Structure d'espace affine sur $\mathcal F$}

\paragraph{2.1\qu  Sous-groupes d'ordre $21$} Soit $M$ un sous-groupe d'ordre 
$7$. Formons le \textit{normalisateur} $\mathrm N(M)$, plus grand 
sous-groupe dans lequel $M$ est distingué.

\textit{Le groupe $M$ coïncide avec son commutant.} En effet, si $h$ commute
avec $f\in \mathcal P\cap M$, comme le polynôme minimal de $f$ est $P(X)$,
de degré $3$, donc égal au polynôme caractéristique, $E$ est cyclique 
pour $f$ et $h$  
appartient à l'algèbre $\mathbf F_2[f]\sim \mathbf F_2[X]/(P)\sim 
\mathbf F_8$
 des applications linéaires 
polynomiales en $f$. Cette algèbre  se déduit 
de  $M$ par adjonction de l'application nulle.

Soit $g\in \mathrm N(M)$. L'automorphisme $h\mapsto ghg^{-1}$ laisse 
$\mathcal P\cap M=\bigl\{f,f^2,f^4\bigr \}$ stable. 

 Si cet automorphisme laisse fixe
 un élément de $\bigl\{f,f^2,f^4\bigr \}$, on peut supposer 
que c'est $f$. On a alors $gf=fg$, donc $g\in M$.

 Si $g\notin M$, cet automorphisme induit une permutation 
circulaire sur $\bigl\{f,f^2,f^4\bigr \}$ et $g^3\in M$. Comme $g$ opère sur 
$\mathrm P(E)$ de cardinal $7$, en le décomposant en cycles disjoints, on 
voit que son ordre ne peut  être $21$, donc est 
nécessairement $3$. Deux cas sont à envisager:

\tir $gfg^{-1}=f^2\qu,\qu gf^2g^{-1}=f^4\qu,\qu gf^4g^{-1}=f$, 

\tir $gfg^{-1}=f^4\qu,\qu gf^2g^{-1}=f\qu,\qu gf^4g^{-1}=f^2$.

\noindent\textbf{Définition} On dira que $g$ est \textit{relié} à $M$ dans le premier cas.
 Il s'agit bien d'une relation entre $g$ 
et $M$ car si elle est vraie pour un $f\in \mathcal P\cap M$, elle est 
aussi vraie pour les autres éléments $f^2,f^4$ de $\mathcal P\cap M$. 

Supposons $g$ relié à $M$. 
Soit $ \varepsilon $ le point fixe 
de $g$. Les relations 
$g( \varepsilon )= \varepsilon $, $gf=f^2g$, 
$gf^2=f^4g$, $gf^4=fg$ impliquent que
 $gf( \varepsilon )=f^2( \varepsilon )$, 
 $gf^2( \varepsilon )=f^4( \varepsilon )$,
 $gf^4( \varepsilon )=f( \varepsilon )$. Autrement dit, la droite stable 
 par $g$ est $\bigl\{f(\varepsilon ),f^2(\varepsilon ),f^4(\varepsilon 
 )\bigr \}$.
 
 Inversement,  soit un point $\varepsilon \in \mathrm P(E)$ et la droite
 $\bigl\{f(\varepsilon ),f^2(\varepsilon ),f^4(\varepsilon )\bigr\}$. 
 Définissons $g$ par $g(\varepsilon )=\varepsilon $, 
 $g\bigl(f(\varepsilon )\bigr)=f^2(\varepsilon )$, 
 $g\bigl(f^2(\varepsilon )\bigr)=f^4(\varepsilon )$. Comme $g$ respecte 
 l'alignement, on a $g\bigl(f^4(\varepsilon )\bigr)=f(\varepsilon )$, 
 $g\bigl(f^3(\varepsilon )\bigr)=f^6(\varepsilon )$, 
 $g\bigl(f^6(\varepsilon )\bigr)=f^5(\varepsilon )$, 
 $g\bigl(f^5(\varepsilon )\bigr)=f^3(\varepsilon )$. D'où l'on déduit que 
 $g\in \mathrm N(M)$ est relié à $M$. On a donc autant d'éléments $g$ d'ordre 
 $3$ reliés à $M$ que de points $\varepsilon \in \mathrm P(E)$, 
 c'est-à-dire $7$. Leurs carrés $g^2$ sont les $7$ autres éléments 
 d'ordre $3$ de $\mathrm N(M)$.
 
 Remarquons que la donnée d'un groupe $M$ d'ordre $7$ détermine sur $\mathrm P(E)$ une 
 structure de $\mathbf F_7$-droite affine pour laquelle $M$ est le groupe 
 des translations. Le groupe $\mathrm N(M)$, de cardinal $21$, apparait alors comme groupe 
 des   homothéties-translations de rapport $1,2,4$ (les carrés de 
 $\mathbf F_7$). C'est le produit semi-direct
 de $M$ par chacun de ses $7$ sous-groupes d'ordre $3$ (Cf. 1.7).
 
\textit{Un élément $g$ d'ordre $3$ est relié à un unique sous-groupe d'ordre 
$7$.} En effet, soit $ \varepsilon $ le point fixe de $g$,
 $\bigl\{ e,g( e),g^2( e)\bigr \}$ la droite laissée stable par $g$,
  alors $g$ est relié au groupe $M=\langle f\rangle$ où $f( \varepsilon 
)= e$, $f^2( \varepsilon )=g( e)$,
 $f^4( \varepsilon )=g^2( e)$. 
Remplacer $f$ 
par $f^2$ ou $f^4$ revient à faire une permutation circulaire sur 
$ e,g( e),g^2( e)$. Ceci se résume dans la définition et la 
proposition suivantes

\noindent\textbf{Définition}\quad \textit{Soit $M$ un groupe d'ordre $7$, 
$\mathcal P\cap M=\bigl\{f,f^2,f^4\bigr \}$. Un élément $g$ d'ordre $3$ est 
\textit{relié} à $M$ s'il vérifie les trois conditions (équivalentes) 
suivantes:}
 $$\Bigl(gfg^{-1}=f^2\Bigr)\Longleftrightarrow 
 \Bigl(gf^2g^{-1}=f^4\Bigr)\Longleftrightarrow \Bigl(gf^4g^{-1}= f\Bigr)$$

\begin{prop} Tout élément $g$ d'ordre $3$ est relié à un unique groupe 
$M$ d'ordre $7$. Si $\varepsilon\in P$ est le point fixe de $g$ et si 
$ \mathcal P\cap M=\bigl\{f,f^2,f^4\bigr \}$, alors 
$\bigl\{f(\varepsilon ),f^2(\varepsilon ),f^4(\varepsilon )\bigr \}$ est la droite 
stable par $g$.

Le groupe $\langle g\rangle =\langle g^2\rangle $ d'ordre $3$, de point 
fixe $\varepsilon $, est relié à 
deux groupes $M$ et $M'$ d'ordre $7$ où, avec $ \mathcal P\cap 
M=\bigl\{f,f^2,f^4\bigr \}$ et $\mathcal P\cap M'=
 \bigl\{f',f'^2,f'^4\bigr \}$, on a 
\begin{eqnarray}&&f_1( \varepsilon )= e\qu , \qu 
f_1^2( \varepsilon )=g( e)\qu ,\qu 
 f_1^4( \varepsilon )=g^2( e)\nonumber\\ 
&&f_2( \varepsilon )= e\qu ,\qu
 f_2^2( \varepsilon )=g^2( e)\qu ,\qu 
 f_2^4( \varepsilon )=g( e)\nonumber \end{eqnarray}
 \end{prop} 

\paragraph{2.2 \qu Les $28$ $3$-sylow(s) et les $28$ paires de $7$-sylow(s)}
 À tout élément $g\in G$ d'ordre $3$, on associe le point fixe 
 $\varepsilon $ et la droite stable orientée $\bigl 
 (e,g(e),g^2(e)\bigr)$, puis les trois éléments $f_1,f_1^2,f_1^4$ 
 d'ordre $7$ appartenant à 
 $\mathcal P$ définis par 
 $$f_1(\varepsilon )=e\qu,\qu f^2_1(\varepsilon )=g(e)\qu,\qu 
 f^4_1(\varepsilon )= g^2(e)$$
On  a ainsi une application $g\mapsto M(g)$ 
de l'ensemble des $56$ éléments d'ordre $3$ sur 
l'ensemble $\mathcal F$ des $8$ sous-groupes d'ordre $7$. 

\begin{prop}$(i)$ L'application de l'ensemble des $28$ sous-groupes d'ordre 
$3$ dans l'ensemble des paires de sous-groupes distincts d'ordre $7$
$$\bigl\{\mathrm I,g,g^2\bigr \}\mapsto \bigl\{M(g),M(g^2)\bigr\}$$
est bijective.

$(ii)$ Étant donné deux couples $(M_1,M_2)$ et $(M'_1,M'_2)$ d'éléments de 
$\mathcal F $ où $M_1\ne M_2$ et $M'_1\ne M'_2$, l'ensemble des $h\in G$ 
tels que $hM_ih^{-1}=M'_i$ pour $i\in \{1,2\}$ est de cardinal $3$.\end{prop}

$(i)$ C'est
une application de l'ensemble des $28$ sous-groupes d'ordre $3$ dans 
l'ensemble des $\mathrm C^2_8= 28$ paires $\bigl\{M_1,M_2\bigr \}$ de sous-groupes 
distincts de $\mathcal F$: à $\Gamma= \bigl\{\mathrm I,g,g^2\bigr \}$ on associe la 
paire $ \bigl\{M(g),M(g^2)\bigr \}$, avec  
$M(g)=\langle f_1\rangle$ et $M(g^2)=\langle f_2\rangle$ où, $\varepsilon 
$ étant  le point fixe de $g$ et $g^2$, 
\begin{eqnarray}&&f_1(\varepsilon )=e\qu ,\qu f_1^2(\varepsilon )=g(e)
\qu ,\qu f_1^4(\varepsilon 
)=g^2(e)\nonumber\\
&&f_2(\varepsilon )=e\qu ,\qu f_2^2(\varepsilon )=g^2(e)\qu ,
\qu f_2^4(\varepsilon )=g(e)\nonumber\end{eqnarray}

Pour montrer que cette application est bijective, 
il suffit de voir qu'elle est injective. 
Pour $g$ d'ordre $3$, les groupes  $M(g)=M_1$ et $M(g^2)=M_2$
sont distincts et $\Gamma =\bigl\{\mathrm I,g,g^2\bigr \}$ est 
dans l'intersection des normalisateurs $\mathrm N(M_1)\cap \mathrm 
N(M_2)$. S'agissant de groupes de cardinal $21$, cette intersection ne 
peut être que de cardinal $3$, donc réduite à $\Gamma$ qui est donc 
l'unique antécédent de la paire $\bigl\{M_1,M_2\bigr \}$.

$(ii)$ On forme $\mathrm N(M_1)\cap \mathrm N(M_2)=\bigl\{\mathrm 
I,g,g^2\bigr \}$ et 
$\mathrm N(M'_1)\cap \mathrm N(M'_2)=\bigl\{\mathrm I,g',g'^2\bigr \}$
 où $g$ (resp. 
$g^2$) est relié à $M_1$ (resp. $M_2$) et $g'$ (resp. $g'^2$) est relié à 
$M'_1$ (resp. $M'_2$). Comme $M_1$ (resp. $M'_1$) est l'unique groupe de 
$\mathcal F $ relié à $g$ (resp. $g'$), l'ensemble des $h$ cherché est 
exactement l'ensemble 
$$\bigl\{h\in G \mid g'=hgh^{-1}\bigr \}$$
Il est bien de cardinal $3$ car le commutant d'un élément $g\in G$ d'ordre 
$3$ est $\langle g\rangle =\bigl\{\mathrm I,g,g^2\bigr \}$ de cardinal $3$.

En particulier, \textit{l'ensemble des $s\in G$ tels que $M_2=sM_1s^{-1}$ et 
$M_1=sM_2s^{-1}$ est constitué des trois involutions du normalisateur 
$\mathrm N(\Gamma)\sim \mathbf S_3$ où $\Gamma =\langle g\rangle = \mathrm N(M_1)\cap 
\mathrm N(M_2)$.}

\paragraph{2.3 \qu Structure de $\mathbf F_2$-espace affine} 
Soit $M_1$ et $M_2$ deux groupes distincts de $\mathcal F$. L'intersection des 
normalisateurs $\mathrm N(M_1)\cap \mathrm N(M_2)$ est un groupe $\Gamma$ 
d'ordre $3$. Si $\varepsilon $ est le point laissé fixe par $\Gamma$, on 
pose $\overrightarrow {M_1M}_2=\overrightarrow {M_2M}_1=\overrightarrow \varepsilon 
\in E$.

 Il s'agit de montrer la \textit{relation de Chasles: si $M,M_1,M_2$ 
sont trois points distincts de $\mathcal F$, alors} 
$$\overrightarrow 
{M_1M_2}=\overrightarrow {M_1M}+\overrightarrow 
{MM_2}$$
Il existe $g_1$ et $g_2$ d'ordre $3$ tels que $g_1$ et $g_2$ soient tous 
deux reliés à 
$M$, $g_1^2$ relié à $M_1$ et $g_2^2$ relié à $M_2$. Autrement dit, 
$$\Gamma _1=\mathrm N(M)\cap \mathrm N(M_1)=\bigl\{\mathrm 
I,g_1,g_1^2\bigr \}\qu,\qu 
\Gamma _2=\mathrm N(M)\cap \mathrm N(M_2)=\bigl\{\mathrm I,g_2,g_2^2\bigr \}$$
pour $f\in \mathcal P\cap M$, $f_1\in \mathcal P\cap M_1$, $f_2\in \mathcal 
P\cap M_2$, on a
$$g_1fg_1^{-1}=g_2fg_2^{-1}=f^2\qu,\qu  g_1^2f_1g_1^{-2}=f_1^2\qu,\qu  
g_2^2f_2g_2^{-2}=f_2^2$$
Si $\varepsilon _1$ et $\varepsilon _2$ sont les points fixes de $g_1$ et 
$g_2$, on a $\overrightarrow {MM_1}=\overrightarrow {\varepsilon _1}$ et 
$\overrightarrow {MM_2}=\overrightarrow {\varepsilon _2}$.

On va montrer qu'il existe deux couples d'involutions $(s'_1,s'_2)$ 
et $(s''_1,s''_2)$ dans $\mathrm N(\Gamma _1)\times \mathrm N(\Gamma _2)$ 
tels que $s'=s'_1s'_2=s'_2s'_1$ et $s''=s''_1s''_2=s''_2s''_1$ soient des 
involutions. On aura alors $M_2= s'M_1s'=s''M_1s''$. On vérifiera que 
$g=s's''$ est d'ordre $3$ appartenant à $\mathrm N(M_1)\cap \mathrm 
N(M_2)$. Si $\varepsilon $ est le point fixe de $g$, on aura 
$\overrightarrow {M_1M_2}=\overrightarrow \varepsilon $. On vérifiera alors la 
relation de Chasles:
 $\overrightarrow {\varepsilon _1}+\overrightarrow {\varepsilon _2}
 +\overrightarrow 
\varepsilon =\overrightarrow 0$.\V

Soit $e\in \mathrm P(E)$. Il existe $i$ et $j$ tels que $\varepsilon _1= 
f^i(e)$ et $\varepsilon _2= f^j(e)$. Les droites stables de $g_1$ et 
$g_2$ sont respectivement $\bigl\{f^{i+1}(e),f^{i+2}(e),f^{i+4}(e)\bigr \}$ et 
$\bigl\{f^{j+1}(e),f^{j+2}(e),f^{j+4}(e)\bigr \}$.
 Ces droites ont un point commun. 
On peut choisir le générateur $f\in \mathcal P\cap M$ et les indices $1$ 
et $2$ de sorte que ce point commun soit 
$f^{i+2}(e)=f^{j+1}(e)$, i.e. que $j=i+1$.
 Pour alléger les notations, on pose $e_k= f^k(e)$. Alors $g_1$ et $g_2$ 
 ont les décompositions en cycles disjoints
 $$g_1=(e_i)(e_{i+1},e_{i+2},e_{i+4})(e_{i+3},e_{i+6},e_{i+5})$$
$$g_2=(e_{i+1})(e_{i+2},e_{i+3},e_{i+5})(e_{i+4},e_{i},e_{i+6})$$
On trouve alors dans les normalisateurs de $\Gamma _1$ et $\Gamma _2$ les 
involutions, chacune composée de deux transpositions disjointes
$$s'_1=(e_{i+2},e_{i+4})(e_{i+6},e_{i+5})\qu,\qu  
s'_2=(e_{i+2},e_{i+5})(e_{i+4},e_{i+6})$$ 
$$s''_1=(e_{i+1},e_{i+4})(e_{i+3},e_{i+5})\qu,\qu  
s''_2=(e_{i+3},e_{i+5})(e_{i},e_{i+6})$$
On a alors
$$s'=s'_1s'_2=s'_2s'_1=(e_{i+2},e_{i+6})(e_{i+4},e_{i+5})$$ 
$$s''=s''_1s''_2=s'_2s''_1=(e_i,e_{i+6})(e_{i+4},e_{i+1}) $$
On en déduit l'élément d'ordre $3$
$$g=s's''=(e_{i+3})(e_i,e_{i+2},e_{i+6})(e_{i+1},e_{i+5},e_{i+4})$$ 
On a donc $\overrightarrow \varepsilon =\overrightarrow {M_1M_2}=
 \overrightarrow {e_{i+3}}$. D'où 
$$\overrightarrow {\varepsilon _1}+\overrightarrow {\varepsilon _2}+
\overrightarrow \varepsilon = 
\overrightarrow {e_i}+ \overrightarrow {e_{i+1}}+
\overrightarrow {e_{i+3}}=\overrightarrow 0$$

Vérifions  que
pour tout $h \in G$, l'application $\widetilde h \colon 
M\mapsto h Mh 
^{-1}$ est affine, d'application linéaire associée $h $. Soit $M_1$ 
et $M_2$ dans $\mathcal  F$, $f_1\in \mathcal P\cap M_1$, $f_2\in 
\mathcal P\cap M_2$, $g$ d'ordre $3$ tels que $gf_1g^{-1}= f_1^2$ et 
$g^2f_2g^{-2}=f_2^2$. Soit $\widetilde h$ l'automorphisme intérieur de $G$ déduit 
de $h $. On aura
 $\widetilde h(g)\widetilde h(f_1)\widetilde h(g)^{-1}
 =\widetilde h(f_1)^2$, 
$\widetilde h(g)^2\widetilde h(f_2)\widetilde h(g)^{-2}
=\widetilde h(f_2)^2$. Si $\varepsilon $ est le 
point fixe de $g$, $h (\varepsilon )$ est le point fixe de 
$\widetilde h(g)$, donc posant $M'_1= \widetilde h(M_1)$
 et $M'_2=\widetilde h(M_2)$, on a bien
$$\overrightarrow {M'_1M'_2}=h (\overrightarrow \varepsilon)=h 
(\overrightarrow {M_1M_2})$$ 

\paragraph{Remarque}  Le dual 
$E^*$ opère canoniquement sur $\mathcal F$ ainsi: étant donné $M_1,M_2$ 
de $\mathcal F$, on forme $\mathrm N(M_1)\cap \mathrm N(M_2)=\bigl\{\mathrm 
I,g,g^2\bigr \}$. Au lieu de considérer le vecteur fixe $\overrightarrow 
\varepsilon $ de $g$, on considère le plan vectoriel stable et la forme 
linéaire $\overleftarrow u \in E^*$ dont il est noyau. On pose alors 
$\overleftarrow {M_1M_2}= \overleftarrow u \in E^*$. On obtient 
une structure de $\mathbf F_2$-espace affine \textquotedblleft{duale}\textquotedblright\ de la 
première. \textit{On n'a pas les mêmes parallélogrammes.}

Soit $M_1,M_2,M'_1,M'_2$ distincts dans $\mathcal F $, $\langle g\rangle =
 \mathrm N(M_1)\cap \mathrm N(M_2)$, $\langle g'\rangle= \mathrm N(M'_1)\cap 
\mathrm N(M'_2)$. Si on avait simultanément
$$\overrightarrow {M_1M_1}=\overrightarrow {M'_1M'_2}=\overrightarrow 
\varepsilon \qu\mbox{ et }\qu 
\overleftarrow{M_1M_2}=\overleftarrow{M'_1M'_2}=\overleftarrow u $$
on aurait $g(\overrightarrow \varepsilon )=g'(\overrightarrow \varepsilon 
)=\overrightarrow \varepsilon $ et $g(\ker u )=g'(\ker u )=\ker 
u $. Comme $E=\Vect(\overrightarrow \varepsilon )\oplus \ker u 
$,  on aurait $g'=g$ ou $g'=g^2$.
Par la proposition 3, on aurait $\bigl\{M_1,M_2\bigr 
\}=\bigl\{M'_1,M'_2\bigr \}$. 

\section{Groupe affine en dimension $3$}

Soit $\mathcal E$ un $\mathbf F_2$-espace affine de dimension $3$, $E$ 
le $\mathbf F_2$-espace vectoriel associé, $\mathcal G$ le groupe affine, 
$G$ le groupe linéaire. Le sous-groupe $\mathcal T$ des 
translations est formé de l'identité et de $7$ translations non nulles 
qui sont involutives.
Les $\mathbf 
F_2$-espaces affines $\mathcal E$ et $\mathcal F$ ont le même espace 
vectoriel associé $E$. En revanche, le groupe affine $\mathcal G$ n'opère 
pas canoniquement sur $\mathcal F$ et il n'existe pas d'isomorphisme 
canonique reliant $\mathcal E$ et $\mathcal F$.

\paragraph{3.1\qu L'espace affine $
\mathcal E'$}
\textit{L'ensemble $\mathcal E'$ des 
isomorphismes affines $\mathcal E'\rightarrow \mathcal F$ induisant 
l'identité sur l'espace $E$ des vecteurs est un espace affine de 
dimension $3$ sur $\mathbf F_2$ admettant $E$ pour espace des vecteurs et 
$\mathcal G$ pour groupe affine.}

Deux isomorphismes $\theta _1,\theta _2\colon \mathcal E\rightarrow 
\mathcal F$ se déduisent l'un de l'autre par une translation de vecteur 
$\overrightarrow v\in E$: pour tout $m\in \mathcal E$, on a 
$\overrightarrow {\theta _1(m)\theta _2(m)}=\overrightarrow v$. On pose 
alors $\overrightarrow {\theta _1\theta _2}=\overrightarrow v$. 

Le groupe $\mathcal G$ opère canoniquement sur $\mathcal E'$: Soit 
$\varphi \in \mathcal G$, $\overline \varphi\in G $ l'application 
linéaire associée et $\theta\in \mathcal E'$. On note encore $\overline 
\varphi $ l'application affine inversible 
 $M\mapsto \overline \varphi M\overline \varphi 
^{-1}$ de $\mathcal F$ sur lui-même. On pose alors $\varphi (\theta )= 
\overline \varphi \circ \theta \circ \varphi ^{-1}$.

Dans ces 
conditions, $\mathcal E$  s'identifie canoniquement à
 \textit{l'espace des isomorphismes
affines $\mathcal E'\rightarrow \mathcal F$ induisant l'identité sur 
$E$}: à $m\in \mathcal E$, on associe l'isomorphisme $\theta \mapsto 
\theta (m)$ de $\mathcal E'$ sur $\mathcal F$. 

\paragraph{3.2\qu Sous-groupes d'ordre $7$ et $21$} Soit $f\in \mathcal G$ tel 
que $\overline f$ soit d'ordre $7$. Alors $1$ n'est pas valeur propre 
donc $f$ a un unique point fixe $m \in \mathcal E$. La surjection 
canonique $\mathcal G \rightarrow G$ induit un isomorphisme $\langle 
f\rangle \sim \langle \overline f\rangle$. Ainsi, l'ensemble $\mathcal 
F\times \mathcal E$ s'identifie à l'ensemble $\mathcal Q$ des 
sous-groupes d'ordre $7$ de $\mathcal G$ en posant pour tout $(M,m)\in 
\mathcal F\times \mathcal E$ 
$$[M,m ]=\bigl\{h\in \mathcal G\mid h(m )=m \mbox{ et } \overline h\in 
M\bigr\}$$

Soit $f$ un générateur de $[M,m ]$. Pour tout $g$ du normalisateur 
$\mathrm N[M,m ]$, $gfg^{-1}\in \mathrm N[M,m ]$. Comme 
$gfg^{-1}$ admet $g(m )$ pour point fixe, $m =g(m )$ car 
$gfg^{-1}$, générateur de $[M,m] $, admet $m $ pour unique point 
fixe. L'application $g\mapsto \overline g$ induit donc un isomorphisme 
entre les normalisateurs $\mathrm N[M,m ]\rightarrow \mathrm N(M)$. 
On a donc 
$$\mathrm N[M,m ]=\bigl\{g\in \mathcal G\mid g(m )=m \mbox{ et } 
\overline g\in \mathrm N(M)\bigr\}$$

\paragraph{3.3\qu Sections du groupe linéaire vers le groupe affine} Soit 
$A,B$ deux groupes de $\mathcal F$, $a,b$ deux points de $\mathcal E$. Si 
les couples $(A,a)$ et $(B,b)$ sont distincts, les groupes $[A,a]$ et 
$[B,b]$ 
sont distincts. L'application $g\mapsto \overline g $ induit des 
isomorphismes $\mathrm N[A,a]\rightarrow \mathrm N(A)$ et $\mathrm 
N[B,b]\rightarrow \mathrm N(B)$. Par ce qui précède, on a
$$\mathrm N[A,a]\cap \mathrm N[B,b]=\bigl\{g\in \mathcal G\mid g(a)=a\qu,\qu 
g(b)=b\qu,\qu \overline g\in \mathrm N(A)\cap \mathrm N(B)\bigr\}$$
Si les couples $(A,a)$ et $(B,b)$ sont distincts, $\mathrm N[A,a]\cap 
\mathrm N[B,b]$ n'est pas réduit à l'identité dans les trois cas suivants:

\tir $A=B$ et $a\ne b$: Soit $\gamma $ l'unique rotation vectorielle 
d'ordre $3$ reliée 
 à $A$  telle que $\gamma (\overrightarrow 
{ab})=\overrightarrow {ab}$ (Cf. 2.1, prop.3). Soit $g\in \mathcal G$ tel 
que $\overline g=\gamma $ et $g(a)=a$. Il est clair que $\langle g\rangle 
= \mathrm N[A,a]\cap \mathrm N[B,b]$.

\tir $A\ne B$ et $a=b$: Soit $\gamma $ tel que $\langle \gamma \rangle 
=\mathrm N (A)\cap \mathrm N(B)$ et $g$ tel que $\overline g=\gamma $ 
et $g(a)=a$. Il est clair que $\langle g\rangle =\mathrm N[A,a]\cap 
\mathrm N[B,b]$.

\tir Supposons $A\ne B$ et $a\ne b$. Pour que $g\in \mathrm N[A,a]\cap 
\mathrm N [B,b]$, il faut et il suffit que $g(a)=a$, $g(b)=b$ et 
$\overline g\in \mathrm N(A)\cap \mathrm N(B)$. On pourra prendre $g$ 
distincte de l'identité si et seulement si $\overrightarrow 
{AB}=\overrightarrow {ab}$ au sens de la structure affine définie sur 
$\mathcal F$ en 2.3.

 Une première classe de sections est formées des  $\mathrm s_a\colon 
f\mapsto f_a$ où $a\in \mathcal E$, $f_a$ étant l'application affine de 
point fixe $a$ et d'application linéaire associée $f$.

Il existe $8$ isomorphismes affines $ \mathcal E\rightarrow 
\mathcal F$ induisant l'identité sur l'espace $E$ des vecteurs. Ils se 
déduisent les uns des autres par des translations. On détermine l'un 
d'eux $\theta $ en fixant l'image $A=\theta (a)
\in \mathcal F$ d'un point $a\in \mathcal E$. 
Pour tout $m\in \mathcal E$, $M=\theta (m)$ est donné par 
$\overrightarrow {AM}=\overrightarrow {am}$.
Considérons la section canonique (Cf.2.3) de $G$ vers le groupe affine de 
$\mathcal F$ associant à $f\in G$ l'application affine
$$\widetilde f\colon \mathcal F\rightarrow \mathcal F\qu,\qu M\mapsto 
fMf^{-1}$$
L'application $\sigma _\theta\colon 
f\mapsto \theta ^{-1}\circ \widetilde f\circ \theta $
 constitue une section appartenant à une autre famille de $8$ 
sections.

Soit $h\in G$ d'ordre $3$ ou $7$. Alors il existe $M\in \mathcal F$ tel 
que $h\in \mathrm N(M)$. On a $hMh^{-1}=M$, donc $\sigma _\theta (h)=h_m$ 
(mais $m=\theta ^{-1}(M)$ n'est pas fixe).

Réciproquement, soit une section $s\colon G\rightarrow \mathcal G$.
Considérons deux sous-groupes $A,B$ distincts de $\mathcal F$. Alors 
$\mathrm N(A)\cap \mathrm N(B)=\langle \gamma \rangle$ où $\gamma $ est 
d'ordre $3$. On a $s(A)=[A,a]$ et $s(B)=[B,b]$ où $(A,a)$ et $(B,b)$ sont 
dans $\mathcal F\times \mathcal E$, ensemble des sous-groupes d'ordre $7$ 
de $\mathcal G$ et $s(\gamma )=g\in \mathrm N[A,a]\cap \mathrm N[B,b]$.
D'après ce qui précède, ou bien $a=b$, ou bien $\overrightarrow 
{AB}=\overrightarrow {ab}$. 

Supposons que $s$ n'est pas  du type $\mathrm s_a$.
 Il existe $A,B$ dans $\mathcal F$ 
tels que $a\ne b$. Soit alors $M\in \mathcal F$ distincts de $A$ et $B$ 
et $s(M)=[M,m]$. Si on avait $m=a$, on aurait $m\ne b$, donc $\overrightarrow 
{BM}=\overrightarrow {bm}=\overrightarrow {ba}=\overrightarrow {BA}$ ce 
qui est impossible car $M\ne A$. Ainsi, $m$ est distinct de $a,b$ et on a
$\overrightarrow {am}=\overrightarrow {AM}$ pour tout $M\in \mathcal F$.

La section $s$ induit donc un isomorphisme 
$\theta \colon \mathcal E\rightarrow \mathcal F$, 
$m\mapsto M$ induisant l'identité sur l'espace $E$ des vecteurs. On voit 
alors que $s=\sigma _\theta $ car $s$ et $\sigma _\theta $ coïncident en les éléments 
d'ordre $3$ ou $7$ qui engendrent le groupe $G$.

Les espaces affines isomorphes $\mathcal E,\mathcal F,\mathcal E'$ ont 
même espace vectoriel $E$, donc même groupe linéaire $G$ et même ensemble 
$\mathcal F$. Les espaces $\mathcal E$ et $\mathcal E'$ jouent le même 
rôle l'un par rapport à l'autre. En revanche, le groupe affine $\mathcal 
G$ de $\mathcal E$ et $\mathcal E'$ ne s'identifie pas canoniquement au 
groupe affine de $\mathcal F$. Par exemple, on a $8$ isomorphismes du groupe 
affine $\mathcal G$ de $\mathcal E$ sur le groupe affine de $\mathcal F$ 
de la forme $f\mapsto \theta f \theta ^{-1}$ où $\theta \in \mathcal E'$.

\paragraph{Remarque} À tout $\theta \in \mathcal E'$ est associée une 
action transitive de $G$ sur $\mathcal E$ via la section $\sigma _\theta 
\colon G\hookrightarrow \mathcal G$. Appliquant la remarque 1.7, on peut 
associer à tout $\theta \in \mathcal E'$ une struture de $\mathbf 
F_7$-droite projective $\Delta _\theta (\mathcal E)$
 sur $\mathcal E$. On a $7$ façons de découper $\mathcal E$ en deux 
 plans $A,A'$ parallèles. On montre que ces $7$ façons $A,A'$ sont les 
 cubes d'une orbite pour $\Delta _\theta (\mathcal E)$ quel que soit 
 $\theta \in \mathcal E'$, l'autre orbite 
 dépendant de $\theta $. Ainsi, 
 l'espace $\mathcal E'$ s'identifie à l'ensemble 
 de ces $8$ structures de $\mathbf F_7$-droites projectives de $\mathcal E$.
 
 \paragraph{3.4 L'ensemble $\mathcal E_1=\mathcal E\cup \mathcal E'$ comme 
$\mathbf F_2$-espace affine de dimension $4$} Soit $\mathcal T_1$ le 
groupe des automorphismes de $\mathcal G$ induisant l'identité sur 
$\mathcal T$ et sur $G=\mathcal G/\mathcal T$. À tout $\overrightarrow 
v\in E$, on associe l'automorphisme intérieur $\widehat {\mathrm 
t}_{\overrightarrow v}\colon \varphi \mapsto \mathrm t_{\overrightarrow 
v}\circ \varphi \circ \mathrm t_{\overrightarrow v}$. On a ainsi une 
inclusion $\mathrm t _{\overrightarrow v}\mapsto 
\widehat {\mathrm t}_{\overrightarrow v}$ de $\mathcal T$ dans $\mathcal T_1$. 

Soit $s\colon G\hookrightarrow \mathcal G$ une section. On en déduit une 
bijection du produit ensembliste $E\times G$ sur $\mathcal G$,   
$ (\overrightarrow v,f)\mapsto \mathrm t _{\overrightarrow 
v}\circ s(f)$.
Pour $\overrightarrow u,\overrightarrow v$ dans $E$ et $f,g$ dans $G$, on a
$$\mathrm t_{\overrightarrow u}\circ s(f)\circ \mathrm t _{\overrightarrow 
v}\circ s(g)=\mathrm t_{\overrightarrow u}\circ
\bigl( s(f)\circ \mathrm t _{\overrightarrow 
v}\circ \circ s(f^{-1})\bigr) \circ s(f)\circ s(g) = \mathrm 
t_{\overrightarrow u} \circ \mathrm t_{f(\overrightarrow v)}\circ s(fg)$$
En transportant la structure de groupe de $\mathcal G$ sur l'ensemble 
$\mathcal T\times G$ par cette bijection, on obtient la loi de groupe
 sur l'ensemble $\mathcal T\times G$ 
$$\bigl (\overrightarrow u,f\bigr )\bigl (\overrightarrow v,g\bigr )=
 \bigl (\overrightarrow 
u+f(\overrightarrow v),fg\bigr )$$
où l'image du sous-groupe $(s(G)$ (resp.  $\mathcal T$) est le 
sous-groupe des $(\overrightarrow 0,f)$ (resp. le sous-groupe distingué 
des $(\overrightarrow v,\mathrm I_E)$).

Un automorphisme $\Phi$ de $ \mathcal G$ induisant 
l'identité sur le sous-groupe $\mathcal T$ et sur le quotient $G=\mathcal 
G/\mathcal T$ se traduit par un automorphisme de $E\times G$ de la forme
$$(\overrightarrow v,f)\mapsto (\overrightarrow v + \overrightarrow 
\varphi (f),f)$$
Le fait qu'il s'agit d'un automorphisme  
se traduit pour l'application 
$\overrightarrow \varphi \colon G\rightarrow E$ par la propriété
$$\forall (f,g)\in G\times G\qu,\qu \overrightarrow \varphi (fg)= 
f\bigl(\overrightarrow \varphi (g)\bigr) + \overrightarrow \varphi (f)$$
L'application $\Phi\mapsto \overrightarrow \varphi $ est une bijection de 
$\mathcal T_1$ sur l'ensemble $E_1$ des applications $G\rightarrow E$ 
vérifiant cette propriété. La structure de $\mathbf F_2$-espace vectoriel 
de $E$ induit une structure de $\mathbf F_2$-espace vectoriel sur $E_1$: 
l'addition de $E_1$ est définie comme suit:
$$\forall (\overrightarrow \varphi ,\overrightarrow \psi)\in E_1\times 
E_1\qu,\qu \forall f\in G\qu,\qu \bigl(\overrightarrow \varphi 
+\overrightarrow \psi\bigr)(f)= \overrightarrow \varphi 
(f)+\overrightarrow \psi (f)$$
On vérifie aussi que la bijection $\Phi\mapsto \overrightarrow \varphi $ 
de $\mathcal T_1$ sur $E_1$ est un isomorphisme de groupes, d'où la 
commutativité de $\mathcal T_1$. L'inclusion $\mathcal T\hookrightarrow 
\mathcal T_ 1$ se traduit par une inclusion $E\hookrightarrow E_1$: au 
vecteur $\overrightarrow v$ est associée l'application $f\mapsto 
\overrightarrow v-f(\overrightarrow v)$ de $G$ vers $E$.

La commutativité de $\mathcal T_1$ permet  de vérifier que cet isomorphisme 
$\mathcal T\rightarrow E_1$, $\Phi \mapsto \overrightarrow \varphi $ est 
\textit{canonique}, i.e. ne 
dépend pas de la section $s\colon G\hookrightarrow \mathcal G$ choisie 
pour le définir.

Le groupe $\mathcal T_1$ opère simplement transitivement sur l'ensemble 
des sections $\mathcal E_1=\mathcal E\cup \mathcal E'$: pour $s\in \mathcal E\cup 
\mathcal E'$ et $\Phi\in \mathcal T_1$, on pose $\Phi(s)= \Phi \circ s$. 

L'ensemble $\mathcal E_1=\mathcal E\cup \mathcal E'$ apparait alors comme 
un $\mathbf F_2$-espace affine de dimension $4$, d'espace vectoriel 
$E_1$, de groupe des translations $\mathcal T_1$, $\mathcal E$ et 
$\mathcal E'$ étant les deux hyperplans parallèles de direction $E$.

Le quotient $E_1/E\sim \mathcal T_1/\mathcal T$, non trivial car de cardinal 
$2$, est le premier groupe de cohomologie $\mathrm H^1_\psi (G,E)$, 
$\psi$ étant l'isomorphisme $G\rightarrow \mathrm{Aut}(E)$ induit par 
l'opération de $G$ sur $E$.

\section{$\mathbf F_2$-espaces affines de dimensions supérieures}

\paragraph{4.1} Soit $\mathcal E$ un $\mathbf F_2$-espace affine de dimension $n\geq 4$, 
de groupe affine $\mathcal G$, de groupe linéaire $G$, d'espace vectoriel 
asscoié $E$. Alors le cardinal de $G$ est égal au nombre de bases, i.e.
\begin{eqnarray}\card G&=&
 (2^n-1)(2^n-2)\cdots (2^n-2^{n-1})\nonumber\\
 &=&
 = (2^n-1)(2^{n-1}-1)\cdots (2^2-1)2^{1+2+\cdots +(n-1)}\nonumber\\
 &=&
 (2^n-1)(2^{n-1}-1)\cdots (2^2-1)2^{\frac{n(n-1)}{2}}\nonumber\\
 \card \mathcal G&=&2^n\card G= 
  (2^n-1)(2^{n-1}-1)\cdots 
  (2^2-1)2^{\frac{n(n+1)}{2}}\nonumber\end{eqnarray}
 
 Supposons qu'\textit{il existe une section $s \colon G\rightarrow \mathcal 
 G$ non du type $\mathrm s_a$}. Montrons alors que $s (G)$ \textit{agit transitivement sur 
 $\mathcal E $.}

 Considérons l'extension $\mathbf F_2\hookrightarrow \mathbf F_{2^n}$. 
 Dans le groupe cyclique $\mathbf F_{2^n}^\times$ existe un élément 
 d'ordre $2^n-1$ de polynôme minimal $P(X)$ irréductible sur $\mathbf 
 F_2$ de degré $n$. Il existe donc dans $G$ des applications d'ordre 
 $2^n-1$ de polynôme minimal $P(X)$. Comme $1$ n'est pas valeur propre 
 de $f$, tout relèvement de $f$ dans $\mathcal G$ a un unique point fixe 
 dans $\mathcal E$. Il existe donc $a\in \mathcal E $ tel que $s(f)=f_a$ et 
 les deux orbites du groupe cyclique $\langle f_a\rangle$ 
  sont $\{a\}$ et $\mathcal E\setminus \{a\}$. 
 Comme $\{a\}$ n'est pas une orbite pour  $s(G)$, l'action de 
 $s(G)$ est bien transitive.

 \paragraph{4.2 Cas  $n=4$}
  On a dans $G$ deux sortes d'applications 
 d'ordre $3$ selon que le polynôme minimal est $X^3-1$ ou $X^2+X+1$. On 
 s'intéresse aux groupes cycliques du deuxième type. On raisonne dans 
 l'espace projectif $P=\mathrm P(E)$ de dimension $3$ et de cardinal 
 $16-1=15$. Un tel groupe $\gamma $ définit $5=\frac{15}{3}$ 
 orbites sur $P$ qui sont des droites disjointes. 
 
 On considère les triplets $\gamma , D_1,D_2$ où $D_1,D_2$ sont deux 
 orbites de $\gamma$. La donnée de $D_1,D_2$ et de deux permutations 
 circulaires sur $D_1,D_2$ détermine $\gamma$. À deux droites disjointes 
 $D_1,D_2$ correspondent deux permutations circulaires sur chaque droite, 
 donc deux groupes ayant $D_1,D_2$ pour orbites.
 
 L'ensemble des droites de $\mathrm P(E)$ est de cardinal 
 $\frac{1}{3}\mathrm C^2_{15}= \frac{14\times 15}{6}=35$. L'ensemble des 
 paires de droites est de cardinal $\mathrm C^2_{35}= 17\times 35$. On a 
 $15$ plans projectifs dans $\mathrm P(E)$ et dans chacun d'eux $7$ 
 droites, donc $\mathrm C^2_7= 21$ paires de droites. On a donc 
 $15\times 21=9\times 35$ paires de droites coplanaires. On a donc 
 $(17-9)\times 35=8\times 35= 280$ paires de droites disjointes.
  D'où $280\times 2=560$ 
 triplets $(\gamma ,D_1,D_2)=(\gamma ,D_2,D_1)$. Pour chaque groupe 
 $\gamma $, on a $5$ orbites, donc $\mathrm C^2_5= 10$ paires $D_1,D_2$.
 On en déduit que le nombre de groupes $\gamma $ est $\frac{560}{10}=56$. 
 
 Comme $1$ n'est pas valeur propre de ces applications d'ordre $3$, tout 
 groupe $\gamma $ se relève en un groupe $\gamma _m$ où $m\in \mathcal E$. 
 S'il existait une section $s\colon G\rightarrow \mathcal G$ sans point 
 fixe, les $16$ stabilisteurs pour l'action de $s(G)$ seraient tous 
 conjugués. Comme $16$ ne divise pas $56$, il est impossible que les 
 groupes $\gamma $ se répartissent équitablement entre les $16$ points de 
 $\mathcal E$.
 
 Ainsi, \textit{le théorème est démontré dans le cas $n=4$.} Ceci permet  
 d'envisager une récurrence.

 \paragraph{4.3 Groupes de Klein de transvections} Soit $\overrightarrow 
 \varepsilon \in E$ et $u\in E^*$, non nuls tels que
  $u(\overrightarrow \varepsilon )=0$. La transvection vectorielle $\mathrm 
  t_{\overrightarrow \varepsilon ,u}$ est définie par
  $$\mathrm t _{\overrightarrow \varepsilon ,u}\colon \overrightarrow 
  x\mapsto \overrightarrow x+u(\overrightarrow x)\overrightarrow 
  \varepsilon $$
  Elle induit l'identité sur $\ker u$ et la translation de vecteur 
  $\overrightarrow \varepsilon $ sur l'hyperplan affine des 
  $\overrightarrow x$ tels que $u(\overrightarrow x)=1$. Elle est 
  involutive.

  On a deux classes de conjugaison de sous-groupes commutatifs
  de transvections de $G$ 
  de cardinal $2^{n-1}$: 
  
  \tir à $u\in E^*$, on associe 
  $\mathrm T_u=\{\mathrm t_{\overrightarrow \varepsilon ,u}\mid 
  \overrightarrow \varepsilon \in E \mbox{ et } u(\overrightarrow 
  \varepsilon )=0\}$
  
  \tir à $\overrightarrow \varepsilon \in E$, on associe 
  $\mathrm T_{\overrightarrow \varepsilon }=\{\mathrm t_{\overrightarrow 
  \varepsilon ,u}\mid u\in E^* \mbox{ et } u(\overrightarrow \varepsilon )=0\}$
  
  \noindent On vérifie que :
 \begin{eqnarray}&&\forall u\in E^*\qu,\qu \forall (\overrightarrow 
 {\varepsilon _1},\overrightarrow {\varepsilon _2})\in E\times E\qu,\qu 
 \mathrm t_{\overrightarrow {\varepsilon _1},u}\circ \mathrm 
 t_{\overrightarrow {\varepsilon _2},u}=\mathrm t_{\overrightarrow 
 {\varepsilon _1}+\overrightarrow {\varepsilon _2},u}\nonumber\\
 &&\forall \overrightarrow \varepsilon \in E\qu,
 \qu \forall (u_1,u_2) \in E^*\times E^*\qu,\qu 
 \mathrm t_{\overrightarrow {\varepsilon },u_1}\circ \mathrm 
 t_{\overrightarrow {\varepsilon },u_2}=\mathrm t_{\overrightarrow 
 \varepsilon,u_1+u_2}\nonumber\end{eqnarray}

  Les applications affines admettant la transvection vectorielle $\mathrm  t   
  _{\overrightarrow \varepsilon ,u}$ pour application linéaire associée
   ne sont en général pas involutives. 
  Celles qui le sont  sont d'un des deux  types suivants:
  
  \tir les $2$ transvections affines du type
   $\mathrm t_{\overrightarrow \varepsilon ,A}$ 
  induisant l'identité sur un hyperplan affine $A$ de direction $\ker u$ 
  et la translation de vecteur $\overrightarrow \varepsilon $ sur 
  l'hyperplan parallèle,
  
  \tir les $2^{n-1}-2$ relèvements affines involutifs sans points fixes
   $\mathrm t_{A,\overrightarrow \alpha ;B,\overrightarrow \beta }$ où 
 $A$ et $B$ sont les hyperplans de direction 
 $\ker u$ et $\overrightarrow \alpha +\overrightarrow 
  \beta =\overrightarrow \varepsilon $, 
  induisant sur les hyperplans affines $A$ et $B$ les 
  translations de vecteurs $\overrightarrow \alpha $ et $\overrightarrow 
  \beta $. 
  
  On a deux types de groupes de Klein de transvections vectorielles: 
  \begin{eqnarray}&&\mathrm K_{\overrightarrow \varepsilon ,u,v,w}=\bigl\{\mathrm I_E, \mathrm t_{\overrightarrow 
  \varepsilon ,u},\mathrm t_{\overrightarrow \varepsilon ,v},\mathrm 
  t_{\overrightarrow \varepsilon ,w}\bigr\} \qu \mbox{ où }\qu 
  u+v+w=0\nonumber\\
  &&\mathrm K_{\overrightarrow \xi ,\overrightarrow \eta ,\overrightarrow 
  \zeta ,u}=\bigl\{\mathrm I_E, \mathrm t_{\overrightarrow 
  \xi  ,u},\mathrm t_{\overrightarrow \eta  ,u},\mathrm 
  t_{\overrightarrow \zeta    ,u}\bigr\} \qu \mbox{ où }\qu 
  \overrightarrow \xi  +\overrightarrow \eta +\overrightarrow \zeta 
  =\overrightarrow 0 \nonumber\end{eqnarray}
 
\paragraph{4.4 Lemme fondamental} \textit{Soit $u,v,w$ des formes linéaires 
 distinctes vérifiant $u+v+w=0$, $\overrightarrow \varepsilon $ un 
 vecteur non nul, $s\colon G\hookrightarrow \mathcal G$ une section,
  $\Gamma=s\bigl (\mathrm K_{\overrightarrow \varepsilon ,u,v,w}\bigr )$.
   Alors, si la dimension est $n\geq 4$, 
  les éléments de $\Gamma$ distincts de l'identité 
  sont des transvections affines 
  d'hyperplans concourants.}\V
  
  Le système de formes linéaires $(u,v,w)$ se relève de $8$ façons en un 
  système de formes affines de somme constante sur $\mathcal E$. On note 
  encore $(u,v,w)$ un des $4$ relèvements de formes affines de somme 
  identiquement nulle sur $\mathcal E$. Considérons les $8$ hyperplans 
  affines
  \begin{eqnarray}U_0=\bigl\{m\in \mathcal E\mid u(m)=0\bigr\}&,&
 U_1=\bigl\{m\in \mathcal E\mid u(m)=1\bigr\}\nonumber\\ 
 V_0=\bigl\{m\in \mathcal E\mid v(m)=0\bigr\}&,&
 V_1=\bigl\{m\in \mathcal E\mid v(m)=1\bigr\}\nonumber\\
  W_0=\bigl\{m\in \mathcal E\mid w(m)=0\bigr\}&,&
 W_1=\bigl\{m\in \mathcal E\mid w(m)=1\bigr\}\nonumber\end{eqnarray}
 Supposons $\Gamma$ formé des relèvements affines
 $$t_u=\mathrm t_{U_0,\overrightarrow {u_0};U_1,\overrightarrow {u_1}}\qu,\qu
 t_v=\mathrm t_{V_0,\overrightarrow {v_0};V_1,\overrightarrow {v_1}}\qu,\qu
t_w= \mathrm t_{W_0,\overrightarrow {w_0};W_1,\overrightarrow {w_1}}$$
\begin{eqnarray}\mbox{avec} &&\overrightarrow {u_0}+\overrightarrow 
{u_1}=\overrightarrow \varepsilon \mbox{ dans } \ker u\nonumber\\
 &&\overrightarrow {v_0}+\overrightarrow 
{v_1}=\overrightarrow \varepsilon \mbox{ dans } \ker v\nonumber\\ 
&&\overrightarrow {w_0}+\overrightarrow 
{w_1}=\overrightarrow \varepsilon \mbox{ dans } \ker 
w\nonumber\end{eqnarray}
On a les deux séries d'égalités
\begin{eqnarray} \overrightarrow {u_0}+\overrightarrow {v_0}+
 \overrightarrow {w_0}=\overrightarrow {u_0}+\overrightarrow 
 {v_1}+\overrightarrow {w_1}=\overrightarrow {u_1}+\overrightarrow 
 {v_0}+\overrightarrow {w_1}=\overrightarrow {u_1}+\overrightarrow 
 {v_1}+\overrightarrow {w_0}=\overrightarrow \alpha \nonumber\\
 \overrightarrow {u_1}+\overrightarrow {v_1}+\overrightarrow {w_1}=
 \overrightarrow {u_1}+\overrightarrow {v_0}+\overrightarrow {w_0}=
 \overrightarrow {u_0}+\overrightarrow {v_1}+\overrightarrow {w_0}=
 \overrightarrow {u_0}+\overrightarrow {v_0}+\overrightarrow 
 {w_1}=\overrightarrow \beta \nonumber\end{eqnarray}
 avec $\overrightarrow \alpha +\overrightarrow \beta =\overrightarrow 
 \varepsilon $.
  L'hyperplan affine 
  $U_0= \bigl(U_0\cap V_0\cap W_0\bigr)\cup \bigl(U_0\cap V_1\cap 
 W_1)$ est stable par $t_u$.
 
 \tir $a)$ Supposons qu'il existe $m\in U_0\cap V_0\cap W_0$ 
 tel que $t_u(m)=m+\overrightarrow {u_0}
 \in U_0\cap V_1\cap W_1$. Alors $\overrightarrow {u_0}\notin \ker v\cup 
 \ker w$ et $v(\overrightarrow {u_0})=w(\overrightarrow {u_0})=1$. On a
 $$m+\overrightarrow {u_0}+\overrightarrow {v_1}=
 \bigl(t_v\circ t_u)(m)=t_w(m)= 
 m+\overrightarrow {w_0}$$
 donc $\overrightarrow \beta =\overrightarrow {u_0}+\overrightarrow 
 {v_1}+\overrightarrow {w_0}=\overrightarrow 0$ et $\overrightarrow 
 \alpha =\overrightarrow \varepsilon $.
 
 \tir $b)$ Supposons qu'il existe $m\in U_0\cap V_0\cap W_0$ 
 tel que $t_u(m)=m+\overrightarrow {u_0}
 \in U_0\cap V_0\cap W_0$. Alors $\overrightarrow {u_0}\in \ker u\cap 
 \ker v\cap \ker w$ et $v(\overrightarrow {u_0})=w(\overrightarrow 
 {u_0})=0$. 
 On a
 $$m+\overrightarrow {u_0}+\overrightarrow {v_0}=
 \bigl(t_v\circ t_u)(m)=t_w(m)= m+\overrightarrow {w_0}$$
 donc $\overrightarrow \alpha  =\overrightarrow {u_0}+\overrightarrow 
 {v_0}+\overrightarrow {w_0}=\overrightarrow 0$ et $\overrightarrow 
 \beta  =\overrightarrow \varepsilon $.
 
 Comme la dimension est $n>3$, $\dim \ker u\geq 3$. Il 
 existe donc une forme linéaire non nulle $v'\ne u$ telle que le plan 
 vectoriel $\{\overrightarrow 0,\overrightarrow {\varepsilon }, 
 \overrightarrow {u_0},\overrightarrow {u_1}\}$ soit contenu dans $\ker 
 v'$, donc aussi dans $\ker w'$ où $w'=u+v'$. Le sous-groupe $\Gamma'=
 s\bigl (\mathrm K_{\overrightarrow \varepsilon ,u,v',w'}\bigr)$ 
 vérifie la condition 
 $b)$. Comme $\mathrm K_{\overrightarrow \varepsilon 
 ,u,v,w}$ 
 est conjugué de $\mathrm K_{\overrightarrow \varepsilon ,u,v',w'}$ dans 
 $G$, $\Gamma'=s\bigl (\mathrm K_{\overrightarrow \varepsilon 
 ,u,v',w'}\bigr)$ et 
 $\Gamma = s\bigl (\mathrm K_{\overrightarrow \varepsilon ,u,v,w}\bigr)$ sont 
 conjugués dans $s(G)$. Ainsi, $\Gamma$ vérifie la condition $b)$. 
 
 Les 
 sous-espaces affines $U_0\cap V_0\cap W_0$ sont donc stables par $t_u$. 
 En faisant varier les formes linéaires $v$ et $w$, on voit finalement 
 que les droites de $U_0$ sont toutes stables par $t_u$, donc que 
 $\overrightarrow {u_0}=\overrightarrow \varepsilon $, $\overrightarrow 
 {u_1}=\overrightarrow 0$, d'où $t_u$ est la transvection affine 
 d'hyperplan $U_0$ de vecteur $\overrightarrow \varepsilon $.
 
 Comme toute transvection vectorielle appartient à un groupe du type 
 $\mathrm K_{\overrightarrow \varepsilon ,u,v,w}$, on déduit du lemme le 
 corollaire ci-dessous:
 
 \noindent\textbf{Corollaire}\quad \textit{L'image d'une transvection 
 vectorielle par une section $s\colon G\hookrightarrow \mathcal G$ est 
 une transvection affine.}\V 
 
 Dans la suite, pour montrer le théorème, on suppose la dimension $n\geq 
 5$, l'hypothèse de récurrence étant que le théorème est supposé vrai en 
 dimension $n-1$. On donne une section $s\colon G\hookrightarrow \mathcal 
 G$. Il s'agit de montrer qu'elle est du type $\mathrm s_a$ où $a\in \mathcal E$.
 
 \paragraph{4.5 Le groupe $\mathrm S_u$ où $u\in E^*$} Soit $u\in E^*$ un 
 forme linéaire non nulle. On forme le sous-groupe de $G$
 $$\mathrm S_u=\{f\in G\mid u\circ f=u\}$$
 Alors $\mathrm S_u$ opère sur les hyperplans stables vectoriel et affine
 $$\mathrm F_{u,0}=\ker u=\{\overrightarrow x\in E\mid u(\overrightarrow 
 x)=0\}\qu,\qu \mathrm F_{u,1}=\{\overrightarrow x\in E\mid 
 u(\overrightarrow x)=1\}$$
 L'action de $\mathrm S_u$ sur $\mathrm F_{u,1}$ permet d'identifier 
 $\mathrm S_u$ au groupe affine. Alors $\mathrm F_{u,0}$ s'identifie à 
 l'espace des vecteurs de $\mathrm F_{u,1}$. L'application $\mathrm 
 S_u\rightarrow \mathrm F_{u,0}$ associant à $f\in \mathrm S_u$ sa 
 restriction à $\mathrm F_{u,0}$ est la surjection canonique du groupe 
 affine sur le groupe linéaire. Le noyau $\mathrm T_u$, groupe des 
 transvections d'hyperplan $\mathrm F_{u,0}$, s'identifie au groupe des 
 translations de l'espace affine $\mathrm F_{u,1}$, la transvection $\mathrm 
 t_{\overrightarrow \varepsilon ,u}$ s'identifiant à la translation de 
 vecteur $\overrightarrow \varepsilon $. 
 
 L'hypothèse de récurrence permet d'affirmer que les sections 
 $\GL(\mathrm F_{u,0})\hookrightarrow  \mathrm S_u$ envoient 
 $\GL(\mathrm F _{u,0})$ en des sous-groupes de $\mathrm S_u$ de la forme
 $$\mathrm G_{\overrightarrow \delta ,u}=\{f\in \mathrm S_u\mid 
 f(\overrightarrow \delta )=\overrightarrow \delta \}=\{f\in G\mid 
 f(\overrightarrow \delta )=\overrightarrow \delta  \mbox{ et } u\circ 
 f=u\}$$
 Soit $A,B$ les deux hyperplans  de $\mathcal E$ de direction $\mathrm 
 F_{u,0}$. 
 Pour qu'une application affine ait son application linéaire associée 
 dans $\mathrm S_u$, il faut et il suffit que les hyperplans $A$ et $B$, 
 ou bien soient stables, ou bien soient échangés. Ces applications 
 affines forment un sous-groupe $\mathcal K_u$ admettant pour sous-groupe 
 d'indice $2$ le sous-groupe $\mathcal H_u$ des applications affines 
 laissant $A$ et $B$ stables.
 
 Montrons que $\mathrm S_u$ n'a pas de sous-groupes d'indice $2$. 
 Comme les groupes linéaires 
 sont simples en dimension au moins $3$, les $\mathrm G_{\overrightarrow 
 \delta ,u}$ sont contenus dans tout sous-groupe d'indice $2$ de $\mathrm 
 S_u$. Il suffit de montrer que  $\mathrm S_u$  est engendré par les sous-groupes $\mathrm 
 G_{\overrightarrow \delta ,u}$.
 
  Soit 
 $\overrightarrow \varepsilon \in \mathrm F _{u,0}$, $\overrightarrow 
 \delta \notin \mathrm F_{u,0}$, $\overrightarrow {\delta _1}= 
 \overrightarrow \delta +\overrightarrow \varepsilon $. Soit $f\in 
 \mathrm G_{\overrightarrow \delta ,u}$ et $f _1\in \mathrm 
 G_{\overrightarrow {\delta _1},u}$ de même application linéaire associée 
 $g$. Alors $f _1\circ f^{-1}$ envoie $\overrightarrow \delta $ en 
 $g(\overrightarrow \delta )+\overrightarrow \varepsilon $ donc est la 
 transvection de vecteur $g(\overrightarrow \delta )+\overrightarrow 
 \delta +\overrightarrow \varepsilon $. Comme 
 $g(\overrightarrow \delta )$ peut être pris quelconque dans $\mathrm 
 F_{u,1}$, toute transvection est dans le sous-groupe engendré par les 
 $\mathrm G_ {\overrightarrow \delta ,u}$. 
 
 La section $s$ envoie donc $\mathrm S_u$ dans le sous-groupe $\mathcal 
 H_u$ des applications affines laissant les hyperplans parallèles $A$ et 
 $B$ stables. Soit $\mathcal A$ et $\mathcal B$ les groupes affines de 
 $A$ et $B$.
 En prenant les restrictions à $A$ et $B$ des $s(f)$ où $f\in \mathrm S 
 _u$, on a des morphismes
  $\mathrm S_u\rightarrow \mathcal A$ et
   $\mathrm S_u\rightarrow \mathcal B$ dont les noyaux ne sont formés que 
   de transvections. Pour tout $\overrightarrow \delta \in \mathrm 
   F_{u,1}$, les restrictions de ces morphismes à $\mathrm 
   G_{\overrightarrow \delta ,u}$ donnent des morphismes injectifs 
   $\mathrm G_{\overrightarrow \delta ,u}\hookrightarrow \mathcal A$ et 
   $\mathrm G_ {\overrightarrow \delta ,u}\hookrightarrow \mathcal B$ 
   qu'on peut considérer comme des sections du groupe linéaire dans le 
   groupe affine en dimension $n-1$. L'hypothèse de récurrence permet 
   d'affirmer qu'il existe $a\in A$ et $b\in B$ fixes par tout $s(g)$ où 
   $g\in \mathrm G_{\overrightarrow \delta ,u}$. 
   
   Ainsi, \textit{pour toute section $s\colon G\hookrightarrow \mathcal G$ et tout 
   sous-groupe $\mathrm G_{\overrightarrow \delta ,u}\subset G$, il existe $a\in 
   A$ et $b\in B$ tels que $\overrightarrow {ab}=\overrightarrow \delta $ 
   et que $s(G_{\overrightarrow \delta ,u})$ soit le sous-groupe $\mathrm 
    G_{a,b,u}\subset \mathcal G$ des applications affines de $\mathcal E$
    laissant fixes   les points $a$ et $b$, 
    stables les hyperplans $A=a+\mathrm F_{u,0}$ et $B=b+\mathrm 
    F_{u,0}$.}
    
    \paragraph{4.6 Lemme} \textit{Soit $s\colon G\hookrightarrow \mathcal 
    G$ une section et $u\in E^*$ une forme linéaire non nulle. Il existe 
    un unique point $a \in \mathcal E$ fixe pour toutes les applications de 
    $s(\mathrm S_u)$.}
 
 Soit $\overrightarrow {\varepsilon _1}$ dans $ \mathrm F_{u,0}$,
  $\overrightarrow 
 \delta $ et $\overrightarrow {\delta _1}=\overrightarrow \delta 
 +\overrightarrow {\varepsilon _1}$ dans $\mathrm F_{u,1}$, on a $\mathrm 
 G_{\overrightarrow {\delta _1},u}= \mathrm t_{\overrightarrow 
 {\varepsilon _1},u}\mathrm G_{\overrightarrow \delta ,u}\mathrm 
 t_{\overrightarrow {\varepsilon _1},u}$.     Par 4.4, $s(\mathrm 
 t_{\overrightarrow \sigma ,u})$ est une transvection de vecteur 
 $\overrightarrow {\varepsilon _1}$ et d'hyperplan $A$ ou $B$. Posant 
 $s(\mathrm G_{\overrightarrow \delta ,u})=\mathrm G_{a,b,u}$ et 
 $s\bigl (\mathrm G_{\overrightarrow {\delta _1},u}\bigr)= \mathrm G_{a_1,b_1,u}$, 
 par conjugaison, la transvection affine
  $s(\mathrm t _{\overrightarrow {\varepsilon _1},u})$ envoie 
 $a$ et $b$ en $a_1$ et $b_1$. On a donc ou bien $a_1=a$ et 
 $b_1=b+\overrightarrow {\varepsilon _1}$, ou bien $a_1=a+\overrightarrow 
 {\varepsilon _1}$ et $b_1=b$.
 
 Soit $\overrightarrow {\varepsilon _2}\ne \overrightarrow {\varepsilon 
 _1}$ dans $ \mathrm F _{u,0}$ et 
 $\overrightarrow {\delta _2}= \overrightarrow {\delta 
 _1}+\overrightarrow {\varepsilon _2}$, $s(\mathrm G_{\overrightarrow 
 {\delta _2},u})= \mathrm G_{a_2,b_2,u}$. Supposons $\overrightarrow 
 {bb_1}=\overrightarrow {\varepsilon _1}$. Alors 
 $\overrightarrow {bb_2}=\overrightarrow {bb_1}+\overrightarrow {b_1b_2}=
 \overrightarrow {\varepsilon _1}+ \overrightarrow {b_1b_2}$
 est égal à $\overrightarrow 0$ ou $\overrightarrow {\varepsilon 
 _1}+\overrightarrow {\varepsilon _2}$. Comme $\overrightarrow {b_1b_2}$ 
 vaut $\overrightarrow 0$ ou $\overrightarrow {\varepsilon _2}$, la seule 
 possibilité est $\overrightarrow {b_1b_2}=\overrightarrow {\varepsilon 
 _2}$ et $\overrightarrow {bb_2}=\overrightarrow {\varepsilon 
 _1}+\overrightarrow {\varepsilon _2}$. On a alors $a=a_1=a_2$. Les 
 groupes $s(\mathrm G_{\overrightarrow \delta ,u})$ sont donc tous de la 
 forme $\mathrm G_{a,b,u}$ où $a\in A$ est indépendant de 
 $\overrightarrow \delta $. Comme les $s(\mathrm G_{\overrightarrow 
 \delta ,u})$ engendrent $s(\mathrm S_u)$, $a$ est bien point fixe de 
 toutes les applications de $s(\mathrm S_u)$.
 
 On obtient l'unicité de $a$ ainsi. Soit $\overrightarrow \delta $ et 
 $\overrightarrow {\delta _1}$ distincts dans $\mathrm F_{u,1}$, $f\in 
 \mathrm G _{\overrightarrow \delta ,u}$ (resp. $f_1\in \mathrm 
 G_{\overrightarrow {\delta _1},u}$) admettant $\overrightarrow \delta $ 
 (resp. $\overrightarrow {\delta _1}$) pour unique vecteur fixe. Alors 
 $s(f)$ (resp. $s(f_1)$) admet $a$ et $b=a+\overrightarrow \delta $ 
 (resp. $a$ et 
 $b_1=a+\overrightarrow {\delta _1}$) pour seuls points fixes. Comme 
 $b\ne b_1$, $a$ est l'unique point fixe commun aux $s(g)$ où $g$ décrit 
$ \mathrm S_u$. Ce point $a$ dépendant à priori de $u$, on le note 
$a_u$.

\paragraph{4.7 Fin de la preuve} Supposons que la section 
$s\colon G\hookrightarrow \mathcal G$ ne soit pas du type $\mathrm s_a$. Par 4.1, 
l'action de $s(G)$ sur $\mathcal E$ serait transitive. Tout point de 
$\mathcal E$ serait d'au moins une façon de la forme $a_u$ où $u\in E^*$, 
$u\ne 0$. C'est impossible car il n'existe que $2^n-1$ formes linéaires 
non nulles et $2^n$ points dans $\mathcal E$.\V\V

 \quad\quad \quad \quad \quad\quad\quad\quad\quad\quad\quad{\Large Bibliographie}
\V\V\
 [1] Arnaudiès  J.M., Bertin J.\qu \textit{Groupes, Algèbres et Géométrie} 
 (Ellipses 1993)
 
 [2] Coxeter H.S.M., \textit{Introduction to geomtry} (J. Wiley and Sons, 
 1989)

 [3] Dieudonné J., \textit{Algèbre linéaire et géométrie élémentaire} (Hermann 1968)

 [4] Hilbert D. et Cohn-Vossen S. \textit{Geometry and the imagination} 
 (Chelsea, 1952)
 
 [5] Mac Lane S. \textit{Homology} (Springer-Verlage 1963)

 [6] Perrin  D., \textit{Cours d'algèbre} (Ellipses, 1996)

 [7] Samuel P., \textit{Géométrie projective} (PUF 1986)

\end{document}